\documentclass[12pt]{article}
\usepackage{amsmath,amsthm,amssymb}
\usepackage{amssymb,latexsym}
\usepackage{mathrsfs}
\newtheorem{theorem}{Theorem}

\newtheorem{remark}{Remark}
\newtheorem{corollary}{Corollary}
\newtheorem{conjecture}{Conjecture}

\begin{document}

\title{\bf  Lower bounds on expressions depending on the functions {\boldmath$\varphi(n)$}, {\boldmath$\psi(n)$} and {\boldmath$\sigma(n)$}, III}

\author{\bf S. I. Dimitrov}

\date{}

\maketitle

\begin{abstract}
This work is concerned with the study of lower bounds for various expressions related to the arithmetic functions $\varphi(n)$, $\psi(n)$ and $\sigma(n)$. Several explicit estimates are established.
\quad\\
{\bf  2020 Math.\ Subject Classification}: 11A25
\end{abstract}

\section{Notations}

Throughout the paper, the symbol $p$ (with or without indices) is reserved for prime numbers.
Let $n \ge 2$ be an integer with canonical prime decomposition
\begin{equation*}
n = \prod_{i=1}^k p_i^{a_i}\,.
\end{equation*}
We write $\Omega(n)$ for the total number of prime factors of $n$, counted with multiplicity, i.e.
\begin{equation*}
\Omega(n) = \sum_{i=1}^k a_i\,, \qquad \Omega(1)=0\,.
\end{equation*}
The Euler totient function is denoted by $\varphi(n)$ and represents the number of integers up to $n$ that are relatively prime to $n$. It satisfies
\begin{equation*}
\varphi(n) = \prod_{i=1}^k p_i^{a_i-1}(p_i - 1)\,, \qquad \varphi(1)=1\,.
\end{equation*}
We also consider the Dedekind psi function, defined multiplicatively by
\begin{equation*}
\psi(n) = \prod_{i=1}^k p_i^{a_i-1}(p_i + 1)\,, \qquad \psi(1)=1.
\end{equation*}
Finally, let $\sigma(n)$ denote the sum-of-divisors function. Its explicit form is given by
\begin{equation*}
\sigma(n) = \prod_{i=1}^k \frac{p_i^{a_i+1}-1}{p_i-1}\,, \qquad \sigma(1)=1\,.
\end{equation*}

\section{Introduction and main results}
\indent

In 2013 Atanassov \cite{Atanassov2013} proved that for every natural number $n\geq2$ the lower bound
\begin{equation*}
\varphi(n)\psi(n)\sigma(n)\geq n^3+n^2-n-1
\end{equation*}
holds. Atanassov's result motivates researchers to investigate lower bounds for expressions involving the functions $\varphi(n)$, $\psi(n)$ and $\sigma(n)$.
More precisely, it has been established that for every integers $n\geq2$, the inequality 
\begin{equation*}
\varphi^k(n)+\psi^k(n)+\sigma^k(n)\geq (n-1)^k +2(n+1)^k
\end{equation*}
holds for $k=2, 3, 4$ in the author's papers \cite{Dimitrov2023, Dimitrov2024a}, and for all $k\geq2$ in S\'{a}ndor and Gryszka \cite{Sandor2025}.
Furthermore, for every integers $n\geq2$, the inequality 
\begin{equation*}
\varphi^k(n)\psi^k(n)+\varphi^k(n)\sigma^k(n)+\sigma^k(n)\psi^k(n)\geq 2(n^2-1)^k+(n+1)^{2k}
\end{equation*}
holds for $k=1, 2$ in \cite{Dimitrov2023, Dimitrov2024a}, and for all $k\geq2$ in \cite{Sandor2025}.
In addition, for every integers $n\geq2$, the inequality 
\begin{align*}
&\varphi^k(n)\big(\psi(n)+\sigma(n)\big)+\psi^k(n)\big(\varphi(n)+\sigma(n)\big)+\sigma^k(n)\big(\varphi(n)+\psi(n)\big)\\
&\geq 2(n+1)[(n-1)^k+2n(n+1)^{k-1}].
\end{align*}
holds for $k=2, 3$ in \cite{Dimitrov2024a}, and for all $k\geq2$ in \cite{Sandor2025}.
Interesting inequalities of this type were also proposed by Mandal \cite{Mandal}. 
A comprehensive overview can be found in the review article \cite{Dimitrov2024b}. 
As a continuation of these studies, we establish lower bounds on new expressions depending on the functions $\varphi(n)$, $\psi(n)$ and $\sigma(n)$.
\begin{theorem}\label{Theorem1}
For every integers $n\geq2$, we have
\begin{equation}\label{lowerbound1}
\big(\varphi(n)+\psi(n)+\sigma(n)\big)\left(\frac{1}{\varphi(n)}+\frac{1}{\psi(n)}+\frac{1}{\sigma(n)}\right)\geq\frac{9n^2-1}{n^2-1}\,.
\end{equation}
\end{theorem}

\begin{theorem}\label{Theorem2}
For every integers $n\geq2$, we have
\begin{equation}\label{lowerbound2}
\big(\varphi(n)+\psi(n)+\sigma(n)\big)\left(\frac{1}{\varphi(n)+\psi(n)}+\frac{1}{\varphi(n)+\sigma(n)}+\frac{1}{\psi(n)+\sigma(n)}\right)\geq\frac{9n^2+9n+2}{2n(n+1)}\,.
\end{equation}
\end{theorem}

\begin{corollary}
As a byproduct of Theorem \ref{Theorem2}, it follows that for every integers $n\geq2$, we have
\begin{equation*}
\frac{\varphi(n)}{\psi(n)+\sigma(n)}+\frac{\psi(n)}{\varphi(n)+\sigma(n)}+\frac{\sigma(n)}{\varphi(n)+\psi(n)}\geq\frac{3n^2+3n+2}{2n(n+1)}\,.
\end{equation*}
\end{corollary}

\begin{theorem}\label{Theorem3}
For every integers $n\geq2$, we have
\begin{align}\label{lowerbound3}
&\Big(\varphi(n)\big(\psi(n)+\sigma(n)\big)+\psi(n)\big(\varphi(n)+\sigma(n)\big)+\sigma(n)\big(\varphi(n)+\psi(n)\big)\Big)\nonumber\\
&\times\left(\frac{1}{\varphi(n)\big(\psi(n)+\sigma(n)\big)}+\frac{1}{\psi(n)\big(\varphi(n)+\sigma(n)\big)}+\frac{1}{\sigma(n)\big(\varphi(n)+\psi(n)\big)}\right)\nonumber\\
&\geq\frac{9n^2-9n+2}{n(n-1)}\,.
\end{align}
\end{theorem}

\begin{theorem}\label{Theorem4}
For every integers $n\geq2$, we have
\begin{equation}\label{lowerbound4}
\frac{\varphi(n)\psi(n)}{\sigma(n)\big(\varphi(n)+\psi(n)\big)}+\frac{\varphi(n)\sigma(n)}{\psi(n)\big(\varphi(n)+\sigma(n)\big)}+\frac{\psi(n)\sigma(n)}{\varphi(n)\big(\psi(n)+\sigma(n)\big)}\geq\frac{3n^2-3n+2}{2n(n-1)}\,.
\end{equation}
\end{theorem}

\begin{theorem}\label{Theorem5}
For every integers $n\geq2$, we have
\begin{equation}\label{lowerbound5}
\frac{\varphi(n)+\psi(n)}{\varphi(n)+2\sigma(n)+\psi(n)}+\frac{\varphi(n)+\sigma(n)}{\varphi(n)+2\psi(n)+\sigma(n)}+\frac{\psi(n)+\sigma(n)}{\psi(n)+2\varphi(n)+\sigma(n)}\geq\frac{6n^2+3n+1}{2n(2n+1)}\,.
\end{equation}
\end{theorem}

\begin{theorem}\label{Theorem6}
For every integers $n\geq2$, we have
\begin{align}\label{lowerbound6}
&\frac{\varphi(n)\big(\psi(n)+\sigma(n)\big)}{\psi(n)\big(\varphi(n)+\sigma(n)\big)+\sigma(n)\big(\varphi(n)+\psi(n)\big)}
+\frac{\psi(n)\big(\varphi(n)+\sigma(n)\big)}{\varphi(n)\big(\psi(n)+\sigma(n)\big)+\sigma(n)\big(\varphi(n)+\psi(n)\big)}\nonumber\\
&+\frac{\sigma(n)\big(\varphi(n)+\psi(n)\big)}{\varphi(n)\big(\psi(n)+\sigma(n)\big)+\psi(n)\big(\varphi(n)+\sigma(n)\big)}\geq\frac{6n^2-3n+1}{2n(2n-1)}\,.
\end{align}
\end{theorem}

In addition, we have the following conjectures.
\begin{conjecture}\label{Conjecture1}
For every integers $k\geq1$ and $n\geq2$, we have
\begin{equation*}
\big(\varphi^k(n)+\psi^k(n)+\sigma^k(n)\big)\left(\frac{1}{\varphi^k(n)}+\frac{1}{\psi^k(n)}+\frac{1}{\sigma^k(n)}\right)\geq\frac{5(n^2-1)^k+2(n+1)^{2k}+2(n-1)^{2k}}{(n^2-1)^k}\,. 
\end{equation*}
\end{conjecture}

\begin{remark} 
Following the method used in the proof of Theorem \ref{Theorem1}, it is not difficult to see that Conjecture \ref{Conjecture1} holds for $k=2$ and $k=3$.
\end{remark}

\begin{conjecture}\label{Conjecture2}
For every integers $k\geq1$ and $n\geq2$, we have
\begin{align*}
&\big((\varphi(n)+\psi(n))^k +(\varphi(n)+\sigma(n))^k+(\psi(n)+\sigma(n))^k\big)\\
&\times\left(\frac{1}{(\varphi(n)+\psi(n))^k}+\frac{1}{(\varphi(n)+\sigma(n))^k}+\frac{1}{(\psi(n)+\sigma(n))^k}\right)\\ 
&\geq\frac{5n^k(n+1)^k+2(n+1)^{2k}+2n^{2k}}{n^k(n+1)^k}\,.
\end{align*}
\end{conjecture}

\begin{conjecture}\label{Conjecture3}
For every integers $k\geq1$ and $n\geq2$, we have
\begin{equation*}
\frac{\varphi^k(n)}{\psi^k(n)+\sigma^k(n)}+\frac{\psi^k(n)}{\varphi^k(n)+\sigma^k(n)}+\frac{\sigma^k(n)}{\varphi^k(n)+\psi^k(n)}\geq\frac{4(n+1)^{2k}+(n^2-1)^k+(n-1)^{2k}}{2\big((n+1)^{2k}+(n^2-1)^k\big)}\,. 
\end{equation*}
\end{conjecture}

\begin{remark} 
Following the method used in the proof of Theorem \ref{Theorem2}, it is not difficult to see that Conjecture \ref{Conjecture3} holds for $k=2$ and $k=3$.
\end{remark}

\begin{conjecture}\label{Conjecture4}
For every integers $k\geq1$ and $n\geq2$, we have
\begin{align*}
&\Big(\varphi^k(n)\big(\psi(n)+\sigma(n)\big)+\psi^k(n)\big(\varphi(n)+\sigma(n)\big)+\sigma^k(n)\big(\varphi(n)+\psi(n)\big)\Big)\\
&\times\left(\frac{1}{\varphi^k(n)\big(\psi(n)+\sigma(n)\big)}+\frac{1}{\psi^k(n)\big(\varphi(n)+\sigma(n)\big)}+\frac{1}{\sigma^k(n)\big(\varphi(n)+\psi(n)\big)}\right)\\
&\geq\frac{5n(n+1)^{k-1}(n-1)^k+2(n-1)^{2k}+2n^2(n+1)^{2k-2}}{n(n+1)^{k-1}(n-1)^k}\,. 
\end{align*}
\end{conjecture}

\begin{conjecture}\label{Conjecture5}
For every integers $k\geq1$ and $n\geq2$, we have
\begin{align*}
&\Big(\varphi^k(n)\big(\psi^k(n)+\sigma^k(n)\big)+\psi^k(n)\big(\varphi^k(n)+\sigma^k(n)\big)+\sigma^k(n)\big(\varphi^k(n)+\psi^k(n)\big)\Big)\\
&\times\left(\frac{1}{\varphi^k(n)\big(\psi^k(n)+\sigma^k(n)\big)}+\frac{1}{\psi^k(n)\big(\varphi^k(n)+\sigma^k(n)\big)}+\frac{1}{\sigma^k(n)\big(\varphi^k(n)+\psi^k(n)\big)}\right)\\
&\geq\frac{10(n-1)^{2k}+7(n^2-1)^k+(n+1)^{2k}}{(n-1)^{2k}+(n^2-1)^k}\,. 
\end{align*}
\end{conjecture}

\begin{conjecture}\label{Conjecture6}
For every integers $k\geq1$ and $n\geq2$, we have
\begin{align*}
&\frac{\varphi^k(n)\psi^k(n)}{\sigma^k(n)\big(\varphi^k(n)+\psi^k(n)\big)}+\frac{\varphi^k(n)\sigma^k(n)}{\psi^k(n)\big(\varphi^k(n)+\sigma^k(n)\big)}+\frac{\psi^k(n)\sigma^k(n)}{\varphi^k(n)\big(\psi^k(n)+\sigma^k(n)\big)}\\
&\geq\frac{4(n-1)^{2k}+(n^2-1)^k+(n+1)^{2k}}{2\big((n-1)^{2k}+(n^2-1)^k\big)}\,. 
\end{align*}
\end{conjecture}

\begin{conjecture}\label{Conjecture7}
For every integers $k\geq1$ and $n\geq2$, we have
\begin{align*}
&\frac{\big(\varphi(n)+\psi(n)\big)^k}{\big(\varphi(n)+\sigma(n)\big)^k+\big(\psi(n)+\sigma(n)\big)^k}+\frac{\big(\varphi(n)+\sigma(n)\big)^k}{\big(\varphi(n)+\psi(n)\big)^k+\big(\psi(n)+\sigma(n)\big)^k}\\
&+\frac{\big(\psi(n)+\sigma(n)\big)^k}{\big(\varphi(n)+\psi(n)\big)^k+\big(\varphi(n)+\sigma(n)\big)^k}\geq\frac{4n^{2k}+(n+1)^{2k}+n^k(n+1)^k}{2\big(n^{2k}+n^k(n+1)^k\big)}\,. 
\end{align*}
\end{conjecture}

\begin{conjecture}\label{Conjecture8}
For every integers $k\geq1$ and $n\geq2$, we have
\begin{align*}
&\frac{\varphi^k(n)\big(\psi(n)+\sigma(n)\big)}{\psi^k(n)\big(\varphi(n)+\sigma(n)\big)+\sigma^k(n)\big(\varphi(n)+\psi(n)\big)}
+\frac{\psi^k(n)\big(\varphi(n)+\sigma(n)\big)}{\varphi^k(n)\big(\psi(n)+\sigma(n)\big)+\sigma^k(n)\big(\varphi(n)+\psi(n)\big)}\\
&+\frac{\sigma^k(n)\big(\varphi(n)+\psi(n)\big)}{\varphi^k(n)\big(\psi(n)+\sigma(n)\big)+\psi^k(n)\big(\varphi(n)+\sigma(n)\big)}\\
&\geq\frac{n(n+1)^{k-1}(n-1)^k+(n-1)^{2k}+4n^2(n+1)^{2k-2}}{2n^2(n+1)^{2k-2}+2n(n+1)^{k-1}(n-1)^k}\,.
\end{align*}
\end{conjecture}

\begin{conjecture}\label{Conjecture9}
For every integers $k\geq1$ and $n\geq2$, we have
\begin{align*}
&\frac{\varphi^k(n)\big(\psi^k(n)+\sigma^k(n)\big)}{\psi^k(n)\big(\varphi^k(n)+\sigma^k(n)\big)+\sigma^k(n)\big(\varphi^k(n)+\psi^k(n)\big)}\\
&+\frac{\psi^k(n)\big(\varphi^k(n)+\sigma^k(n)\big)}{\varphi^k(n)\big(\psi^k(n)+\sigma^k(n)\big)+\sigma^k(n)\big(\varphi^k(n)+\psi^k(n)\big)}\\
&+\frac{\sigma^k(n)\big(\varphi^k(n)+\psi^k(n)\big)}{\varphi^k(n)\big(\psi^k(n)+\sigma^k(n)\big)+\psi^k(n)\big(\varphi^k(n)+\sigma^k(n)\big)}\\
&\geq\frac{5(n-1)^{2k}+5(n^2-1)^k+2(n+1)^{2k}}{\big(3(n-1)^k+(n+1)^k\big)\big((n-1)^k+(n+1)^k\big)}\,. 
\end{align*}
\end{conjecture}

\section{Proof of Theorem \ref{Theorem1}}

Consider several cases.

Case 1. \; $\Omega(n)=1$. Bearing in mind that $n$ is a prime number, we write
\begin{equation*}
\big(\varphi(n)+\psi(n)+\sigma(n)\big)\left(\frac{1}{\varphi(n)}+\frac{1}{\psi(n)}+\frac{1}{\sigma(n)}\right)=(3n+1)\left(\frac{1}{n-1}+\frac{2}{n+1}\right)=\frac{9n^2-1}{n^2-1}\,.
\end{equation*}

Case 2. \; $\Omega(n)=2$, $n=pq$, where $p$ and $q$ are distinct primes. We have
\begin{align*}
&\big(\varphi(n)+\psi(n)+\sigma(n)\big)\left(\frac{1}{\varphi(n)}+\frac{1}{\psi(n)}+\frac{1}{\sigma(n)}\right)\\
&=\big((p-1)(q-1)+2(p+1)(q+1)\big)\left(\frac{1}{(p-1)(q-1)}+\frac{2}{(p+1)(q+1)}\right)\\
&=\frac{9p^2q^2+16pq-p^2-q^2+9}{(p^2 - 1) (q^2 - 1)}>\frac{9p^2q^2-1}{p^2q^2-1}=\frac{9n^2-1}{n^2-1}\,.
\end{align*}

Case 3. \; $\Omega(n)=2$, $n=p^2$, where $p$ is a prime. We write
\begin{align*}
&\big(\varphi(n)+\psi(n)+\sigma(n)\big)\left(\frac{1}{\varphi(n)}+\frac{1}{\psi(n)}+\frac{1}{\sigma(n)}\right)\\
\end{align*}
\begin{align*}
&=\big(p(p-1)+p(p+1)+p^2+p+1\big)\left(\frac{1}{p(p-1)}+\frac{1}{p(p+1)}+\frac{1}{p^2+p+1}\right)\\    
&=\frac{(3p^2+p+1)(3p^2+2p+1)}{(p^2-1) (p^2 + p + 1)}>\frac{9p^4-1}{p^4-1}=\frac{9n^2-1}{n^2-1}\,.
\end{align*}
Now we assume that \eqref{lowerbound1} is true for every positive integer $n$ with $\Omega(n)=m$ for some positive integer $m\geq2$.
Let $p$ be a prime number. Then $\Omega(np)=\Omega(n)+1$.

Case A.  \;  $p\nmid n$. Using that
\begin{equation*}
\varphi(n)\leq n-1\,, \qquad \psi(n)\geq n+1 \,, \qquad \sigma(n)\geq n+1
\end{equation*}
we obtain
\begin{align*}
&\big(\varphi(np)+\psi(np)+\sigma(np)\big)\left(\frac{1}{\varphi(np)}+\frac{1}{\psi(np)}+\frac{1}{\sigma(np)}\right)\\
&=\Big((p+1)\big(\varphi(n)+\psi(n)+\sigma(n)\big)-2\varphi(n)\Big)\\
&\times\left(\frac{1}{p+1}\left(\frac{1}{\varphi(n)}+\frac{1}{\psi(n)}+\frac{1}{\sigma(n)}\right)+\frac{2}{(p^2-1)\varphi(n)}\right)\\
&=\big(\varphi(n)+\psi(n)+\sigma(n)\big)\left(\frac{1}{\varphi(n)}+\frac{1}{\psi(n)}+\frac{1}{\sigma(n)}\right)\\
&+\frac{2}{p-1}\frac{\psi(n)+\sigma(n)}{\varphi(n)}-\frac{2}{p+1}\left(\frac{1}{\psi(n)}+\frac{1}{\sigma(n)}\right)\varphi(n)\\
&\geq\frac{9n^2-1}{n^2-1}+\frac{4\left(n+1\right)}{\left(p-1\right)\left(n-1\right)}-\frac{4\left(n-1\right)}{\left(p+1\right)\left(n+1\right)}>\frac{9n^2p^2-1}{n^2p^2-1}\,.
\end{align*}

Case B.  \;  $p\,|\,n$.  We have
\begin{equation}\label{phipsisigmanp}
\varphi(np)=p\varphi(n)\,, \qquad \psi(np)=p\psi(n)\,, \qquad \sigma(np)>p\sigma(n)
\end{equation}
and
\begin{equation}\label{phi<psi<sigma}
\varphi(n)\leq\psi(n)\,, \qquad \psi(n)\leq \sigma(n)\,.
\end{equation}
It is easy to see that \eqref{phipsisigmanp} and \eqref{phi<psi<sigma} imply
\begin{equation*}
\left(\frac{\sigma(np)}{p}-\sigma(n)\right)\left(\frac{1}{\varphi(n)}+\frac{1}{\psi(n)}\right)+\left(\frac{p}{\sigma(np)}-\frac{1}{\sigma(n)}\right)\big(\varphi(n)+\psi(n)\big)\geq0
\end{equation*}
which, together with \eqref{phipsisigmanp} yields
\begin{align*}
&\big(\varphi(np)+\psi(np)+\sigma(np)\big)\left(\frac{1}{\varphi(np)}+\frac{1}{\psi(np)}+\frac{1}{\sigma(np)}\right)\\
&=\Big[p\big(\varphi(n)+\psi(n)+\sigma(n)\big)+\sigma(np)-p\sigma(n)\Big]\\
&\times\left[\frac{1}{p}\left(\frac{1}{\varphi(n)}+\frac{1}{\psi(n)}+\frac{1}{\sigma(n)}\right)+\frac{1}{\sigma(np)}-\frac{1}{p\sigma(n)}\right]\\
\end{align*}
\begin{align*}
&=\big(\varphi(n)+\psi(n)+\sigma(n)\big)\left(\frac{1}{\varphi(n)}+\frac{1}{\psi(n)}+\frac{1}{\sigma(n)}\right)\\
&+\left(\frac{\sigma(np)}{p}-\sigma(n)\right)\left(\frac{1}{\varphi(n)}+\frac{1}{\psi(n)}\right)+\left(\frac{p}{\sigma(np)}-\frac{1}{\sigma(n)}\right)\big(\varphi(n)+\psi(n)\big)\\
&>\frac{9n^2-1}{n^2-1}>\frac{9n^2p^2-1}{n^2p^2-1}\,.
\end{align*}
This completes the proof of Theorem \ref{Theorem1}.

\section{Proof of Theorem \ref{Theorem2}}

Consider several cases.

Case 1. \; $\Omega(n)=1$. Taking into account that $n$ is a prime number, we  have
\begin{align*}
&\big(\varphi(n)+\psi(n)+\sigma(n)\big)\left(\frac{1}{\varphi(n)+\psi(n)}+\frac{1}{\varphi(n)+\sigma(n)}+\frac{1}{\psi(n)+\sigma(n)}\right)\\
&=(3n+1)\left(\frac{1}{n}+\frac{1}{2(n+1)}\right)=\frac{9n^2+9n+2}{2n(n+1)}\,.
\end{align*}

Case 2. \; $\Omega(n)=2$, $n=pq$, where $p$ and $q$ are distinct primes. Then
\begin{align*}
&\big(\varphi(n)+\psi(n)+\sigma(n)\big)\left(\frac{1}{\varphi(n)+\psi(n)}+\frac{1}{\varphi(n)+\sigma(n)}+\frac{1}{\psi(n)+\sigma(n)}\right)\\
&=\big((p-1)(q-1)+2(p+1)(q+1)\big)\left(\frac{2}{(p-1)(q-1)+(p+1)(q+1)}+\frac{1}{2(p+1)(q+1)}\right)\\   
&=\frac{9p^2q^2+9p^2q+9pq^2+2p^2+22pq+2q^2+9p+9q+9}{2(p+1)(q+1)(pq+1)}\\
&>\frac{9p^2q^2+9pq+2}{2pq(pq+1)}=\frac{9n^2+9n+2}{2n(n+1)}\,.
\end{align*}

Case 3. \; $\Omega(n)=2$, $n=p^2$, where $p$ is a prime. We get
\begin{align*}
&\big(\varphi(n)+\psi(n)+\sigma(n)\big)\left(\frac{1}{\varphi(n)+\psi(n)}+\frac{1}{\varphi(n)+\sigma(n)}+\frac{1}{\psi(n)+\sigma(n)}\right)\\
&=\big(p(p-1)+p(p+1)+p^2+p+1\big)\\
&\times\left(\frac{1}{p(p-1)+p(p+1)}+\frac{1}{p(p-1)+p^2+p+1}+\frac{1}{p(p+1)+p^2+p+1}\right)\\    
&=\frac{36p^6+36p^5+44p^4+22p^3+13p^2+3p+1}{2p^2(2p^2+1) (2p^2 + 2p + 1)}>\frac{9p^4+9p^2+2}{2p^2(p^2+1)}=\frac{9n^2+9n+2}{2n(n+1)}\,.
\end{align*}
Now we assume that \eqref{lowerbound2} is true for every positive integer $n$ with $\Omega(n)=m$ for some positive integer $m\geq2$.
Let $p$ be a prime number. Then $\Omega(np)=\Omega(n)+1$.

Case A.  \;  $p\nmid n$. It follows easily from \eqref{phi<psi<sigma} that
\begin{align*}
&\frac{\varphi(n)\sigma(n)}{\big((p-1)\varphi(n)+(p+1)\psi(n)\big)\big(\varphi(n)+\psi(n)\big)}+\frac{\varphi(n)\psi(n)}{\big((p-1)\varphi(n)+(p+1)\sigma(n)\big)\big(\varphi(n)+\sigma(n)\big)}\\
&-\frac{\varphi(n)}{(p+1)\big(\psi(n)+\sigma(n)\big)}\geq0
\end{align*}
which leads to
\begin{align*}
&\big(\varphi(np)+\psi(np)+\sigma(np)\big)\left(\frac{1}{\varphi(np)+\psi(np)}+\frac{1}{\varphi(np)+\sigma(np)}+\frac{1}{\psi(np)+\sigma(np)}\right)\\
&=\Big((p-1)\varphi(n)+(p+1)\psi(n)+(p+1)\sigma(n)\Big)\\
&\times\left(\frac{1}{(p-1)\varphi(n)+(p+1)\psi(n)}+\frac{1}{(p-1)\varphi(n)+(p+1)\sigma(n)}+\frac{1}{(p+1)\big(\psi(n)+\sigma(n)\big)}\right)\\
&=\Big((p+1)\big(\varphi(n)+\psi(n)+\sigma(n)\big)-2\varphi(n)\Big)\\
&\times\frac{1}{p+1}\left(\frac{1}{\varphi(n)+\psi(n)}+\frac{1}{\varphi(n)+\sigma(n)}+\frac{1}{\psi(n)+\sigma(n)}\right)\\
&+\Big((p-1)\varphi(n)+(p+1)\psi(n)+(p+1)\sigma(n)\Big)\\
&\times\left(\frac{1}{(p-1)\varphi(n)+(p+1)\psi(n)}+\frac{1}{(p-1)\varphi(n)+(p+1)\sigma(n)}\right)\\
&-\Big((p+1)\big(\varphi(n)+\psi(n)+\sigma(n)\big)-2\varphi(n)\Big)\frac{1}{p+1}\left(\frac{1}{\varphi(n)+\psi(n)}+\frac{1}{\varphi(n)+\sigma(n)}\right)\\
&\geq\frac{9n^2+9n+2}{2n(n+1)}+\frac{2\varphi(n)\sigma(n)}{\big((p-1)\varphi(n)+(p+1)\psi(n)\big)\big(\varphi(n)+\psi(n)\big)}\\
&+\frac{2\varphi(n)\psi(n)}{\big((p-1)\varphi(n)+(p+1)\sigma(n)\big)\big(\varphi(n)+\sigma(n)\big)}-\frac{2\varphi(n)}{(p+1)\big(\psi(n)+\sigma(n)\big)}\\
&>\frac{9n^2p^2+9np+2}{2np(np+1)}\,.
\end{align*}

Case B.  \;  $p\,|\,n$. It is easy to see that \eqref{phipsisigmanp} and \eqref{phi<psi<sigma} give us
\begin{align*}
\frac{\sigma(np)-p\sigma(n)}{p\big(\varphi(n)+\psi(n)\big)}-\frac{\big(\sigma(np)-p\sigma(n))\psi(n)}{\big(p\varphi(n)+\sigma(np)\big)\big(\varphi(n)+\sigma(n)\big)}
-\frac{\big(\sigma(np)-p\sigma(n))\varphi(n)}{\big(p\psi(n)+\sigma(np)\big)\big(\sigma(n)+\psi(n)\big)}\geq0
\end{align*}
which, together with \eqref{phipsisigmanp} implies
\begin{align*}
&\big(\varphi(np)+\psi(np)+\sigma(np)\big)\left(\frac{1}{\varphi(np)+\psi(np)}+\frac{1}{\varphi(np)+\sigma(np)}+\frac{1}{\psi(np)+\sigma(np)}\right)\\
\end{align*}
\begin{align*}
&=\big(p\varphi(n)+p\psi(n)+\sigma(np)\big)\left(\frac{1}{p\varphi(n)+p\psi(n)}+\frac{1}{p\varphi(n)+\sigma(np)}+\frac{1}{p\psi(n)+\sigma(np)}\right)\\
&=\Big(p\big(\varphi(n)+\psi(n)+\sigma(n)\big)+\sigma(np)-p\sigma(n)\Big)\\
&\times\frac{1}{p}\left(\frac{1}{\varphi(n)+\psi(n)}+\frac{1}{\varphi(n)+\sigma(n)}+\frac{1}{\psi(n)+\sigma(n)}\right)\\
&+\big(p\varphi(n)+p\psi(n)+\sigma(np)\big)\left(\frac{1}{p\varphi(n)+\sigma(np)}+\frac{1}{p\psi(n)+\sigma(np)}\right)\\
&-\Big(p\big(\varphi(n)+\psi(n)+\sigma(n)\big)+\sigma(np)-p\sigma(n)\Big)\frac{1}{p}\left(\frac{1}{\varphi(n)+\sigma(n)}+\frac{1}{\psi(n)+\sigma(n)}\right)\\
&\geq\frac{9n^2+9n+2}{2n(n+1)}+\frac{\sigma(np)-p\sigma(n)}{p\big(\varphi(n)+\psi(n)\big)}-\frac{\big(\sigma(np)-p\sigma(n))\psi(n)}{\big(p\varphi(n)+\sigma(np)\big)\big(\varphi(n)+\sigma(n)\big)}\\
&-\frac{\big(\sigma(np)-p\sigma(n))\varphi(n)}{\big(p\psi(n)+\sigma(np)\big)\big(\sigma(n)+\psi(n)\big)}\\
&>\frac{9n^2p^2+9np+2}{2np(np+1)}\,.
\end{align*}
This completes the proof of Theorem \ref{Theorem2}.

\section{Proof of Theorem \ref{Theorem3}}

Consider several cases.

Case 1. \; $\Omega(n)=1$. Since $n$ is a prime number, we write
\begin{align*}
&\Big(\varphi(n)\big(\psi(n)+\sigma(n)\big)+\psi(n)\big(\varphi(n)+\sigma(n)\big)+\sigma(n)\big(\varphi(n)+\psi(n)\big)\Big)\\
&\times\left(\frac{1}{\varphi(n)\big(\psi(n)+\sigma(n)\big)}+\frac{1}{\psi(n)\big(\varphi(n)+\sigma(n)\big)}+\frac{1}{\sigma(n)\big(\varphi(n)+\psi(n)\big)}\right)\\
&=\big(2(n^2-1)+4n(n+1)\big)\left(\frac{1}{2(n^2-1)}+\frac{1}{n(n+1)}\right)=\frac{9n^2-9n+2}{n(n-1)}\,.
\end{align*}

Case 2. \; $\Omega(n)=2$, $n=pq$, where $p$ and $q$ are distinct primes. After straightforward though lengthy calculations, we deduce
\begin{align*}
&\Big(\varphi(n)\big(\psi(n)+\sigma(n)\big)+\psi(n)\big(\varphi(n)+\sigma(n)\big)+\sigma(n)\big(\varphi(n)+\psi(n)\big)\Big)\\
&\times\left(\frac{1}{\varphi(n)\big(\psi(n)+\sigma(n)\big)}+\frac{1}{\psi(n)\big(\varphi(n)+\sigma(n)\big)}+\frac{1}{\sigma(n)\big(\varphi(n)+\psi(n)\big)}\right)\\
&=\big(4(p^2-1)(q^2-1)+2(p+1)^2(q+1)^2\big)\left(\frac{1}{2(p^2-1)(q^2-1)}+\frac{1}{(p+1)(q+1)(pq+1)}\right)\\   
\end{align*}
\begin{align*}
&>\frac{9p^2q^2-9pq+2}{pq(pq-1)}=\frac{9n^2-9n+2}{n(n-1)}\,.
\end{align*}

Case 3. \; $\Omega(n)=2$, $n=p^2$, where $p$ is a prime. We derive
\begin{align*}
&\Big(\varphi(n)\big(\psi(n)+\sigma(n)\big)+\psi(n)\big(\varphi(n)+\sigma(n)\big)+\sigma(n)\big(\varphi(n)+\psi(n)\big)\Big)\\
&\times\left(\frac{1}{\varphi(n)\big(\psi(n)+\sigma(n)\big)}+\frac{1}{\psi(n)\big(\varphi(n)+\sigma(n)\big)}+\frac{1}{\sigma(n)\big(\varphi(n)+\psi(n)\big)}\right)\\
&=\big(p(p-1)(p(p+1)+p^2+p+1) +  p(p+1)(p(p-1)+p^2+p+1) \\
&+ (p^2+p+1)(p(p-1)+p(p+1))\big)\\
&\times\left(\frac{1}{p(p-1)(p(p+1)+p^2+p+1)}+\frac{1}{p(p+1)(p(p-1)+p^2+p+1)}\right.\\
&\left.+\frac{1}{(p^2+p+1)(p(p-1)+p(p+1))}\right)\\    
&=\frac{16p^{12}+40p^{11}+48p^{10}+24p^9-4p^8-10p^7+2p^6+6p^5+p^4-5p^3-6p^2-3p-1}{4p^8+8p^7+8p^6+2p^5-5p^4-7p^3-6p^2-3p-1}\\
&>\frac{9p^4-9p^2+2}{p^2(p^2-1)}=\frac{9n^2-9n+2}{n(n-1)}\,.
\end{align*}
Now we assume that \eqref{lowerbound3} is true for every positive integer $n$ with $\Omega(n)=m$ for some positive integer $m\geq2$.
Let $p$ be a prime number. Then $\Omega(np)=\Omega(n)+1$.

Case A.  \;  $p\nmid n$. Using \eqref{phi<psi<sigma} together with lengthy but elementary calculations, we obtain
\begin{align*}
&\frac{(p-1)\varphi(n)\big(\psi(n)+\sigma(n)\big)}{\psi(n)\big((p-1)\varphi(n)+(p+1)\sigma(n)\big)}+\frac{\psi(n)\sigma(n)}{(p-1)\varphi(n)\big(\psi(n)+\sigma(n)\big)}
-\frac{\varphi(n)\big(\psi(n)+\sigma(n)\big)}{\psi(n)\big(\varphi(n)+\sigma(n)\big)}\geq0\,,\\
&\frac{(p-1)\varphi(n)\big(\psi(n)+\sigma(n)\big)}{\sigma(n)\big((p-1)\varphi(n)+(p+1)\psi(n)\big)}+\frac{\psi(n)\sigma(n)}{(p-1)\varphi(n)\big(\psi(n)+\sigma(n)\big)}
-\frac{\varphi(n)\big(\psi(n)+\sigma(n)\big)}{\sigma(n)\big(\varphi(n)+\psi(n)\big)}\geq0\,,\\
&\frac{\psi(n)\big((p-1)\varphi(n)+(p+1)\sigma(n)\big)}{\sigma(n)\big((p-1)\varphi(n)+(p+1)\psi(n)\big)}+\frac{\psi(n)\sigma(n)}{(p-1)\varphi(n)\big(\psi(n)+\sigma(n)\big)}
-\frac{\psi(n)\big(\varphi(n)+\sigma(n)\big)}{\sigma(n)\big(\varphi(n)+\psi(n)\big)}\geq0\,,\\
&\frac{\sigma(n)\big((p-1)\varphi(n)+(p+1)\psi(n)\big)}{\psi(n)\big((p-1)\varphi(n)+(p+1)\sigma(n)\big)}+\frac{\psi(n)\sigma(n)}{(p-1)\varphi(n)\big(\psi(n)+\sigma(n)\big)}
-\frac{\sigma(n)\big(\varphi(n)+\psi(n)\big)}{\psi(n)\big(\varphi(n)+\sigma(n)\big)}\geq0
\end{align*}
which yields
\begin{align*}
&\Big(\varphi(np)\big(\psi(np)+\sigma(np)\big)+\psi(np)\big(\varphi(np)+\sigma(np)\big)+\sigma(np)\big(\varphi(np)+\psi(np)\big)\Big)\\
&\times\left(\frac{1}{\varphi(np)\big(\psi(np)+\sigma(np)\big)}+\frac{1}{\psi(np)\big(\varphi(np)+\sigma(np)\big)}+\frac{1}{\sigma(np)\big(\varphi(np)+\psi(np)\big)}\right)\\
\end{align*}
\begin{align*}
&=\Big((p^2-1)\varphi(n)\big(\psi(n)+\sigma(n)\big)+(p+1)\psi(n)\big((p-1)\varphi(n)+(p+1)\sigma(n)\big)\\
&+(p+1)\sigma(n)\big((p-1)\varphi(n)+(p+1)\psi(n)\big)\Big)\\
&\times\left(\frac{1}{(p^2-1)\varphi(n)\big(\psi(n)+\sigma(n)\big)}+\frac{1}{(p+1)\psi(n)\big((p-1)\varphi(n)+(p+1)\sigma(n)\big)}\right.\\
&\left.+\frac{1}{(p+1)\sigma(n)\big((p-1)\varphi(n)+(p+1)\psi(n)\big)}\right)\\
&=\bigg((p^2-1)\Big(\varphi(n)\big(\psi(n)+\sigma(n)\big)+\psi(n)\big(\varphi(n)+\sigma(n)\big)+\sigma(n)\big(\varphi(n)+\psi(n)\big)\Big)\\
&+4(p+1)\sigma(n)\psi(n)\bigg)\frac{1}{p^2-1}\\
&\times\left(\frac{1}{\varphi(n)\big(\psi(n)+\sigma(n)\big)}+\frac{1}{\psi(n)\big(\varphi(n)+\sigma(n)\big)}+\frac{1}{\sigma(n)\big(\varphi(n)+\psi(n)\big)}\right)\\
&+\Big((p^2-1)\varphi(n)\big(\psi(n)+\sigma(n)\big)+(p+1)\psi(n)\big((p-1)\varphi(n)+(p+1)\sigma(n)\big)\\
&+(p+1)\sigma(n)\big((p-1)\varphi(n)+(p+1)\psi(n)\big)\Big)\\
&\times\left(\frac{1}{(p+1)\psi(n)\big((p-1)\varphi(n)+(p+1)\sigma(n)\big)}+\frac{1}{(p+1)\sigma(n)\big((p-1)\varphi(n)+(p+1)\psi(n)\big)}\right)\\
&-\bigg((p^2-1)\Big(\varphi(n)\big(\psi(n)+\sigma(n)\big)+\psi(n)\big(\varphi(n)+\sigma(n)\big)+\sigma(n)\big(\varphi(n)+\psi(n)\big)\Big)\\
&+4(p+1)\sigma(n)\psi(n)\bigg)\left(\frac{1}{\psi(n)\big(\varphi(n)+\sigma(n)\big)}+\frac{1}{\sigma(n)\big(\varphi(n)+\psi(n)\big)}\right)\frac{1}{p^2-1}\\
&\geq\frac{9n^2-9n+2}{n(n-1)}+\\
&+\frac{(p-1)\varphi(n)\big(\psi(n)+\sigma(n)\big)}{\psi(n)\big((p-1)\varphi(n)+(p+1)\sigma(n)\big)}+\frac{\psi(n)\sigma(n)}{(p-1)\varphi(n)\big(\psi(n)+\sigma(n)\big)}
-\frac{\varphi(n)\big(\psi(n)+\sigma(n)\big)}{\psi(n)\big(\varphi(n)+\sigma(n)\big)}\nonumber\\
&+\frac{(p-1)\varphi(n)\big(\psi(n)+\sigma(n)\big)}{\sigma(n)\big((p-1)\varphi(n)+(p+1)\psi(n)\big)}+\frac{\psi(n)\sigma(n)}{(p-1)\varphi(n)\big(\psi(n)+\sigma(n)\big)}
-\frac{\varphi(n)\big(\psi(n)+\sigma(n)\big)}{\sigma(n)\big(\varphi(n)+\psi(n)\big)}\nonumber\\
&+\frac{\psi(n)\big((p-1)\varphi(n)+(p+1)\sigma(n)\big)}{\sigma(n)\big((p-1)\varphi(n)+(p+1)\psi(n)\big)}+\frac{\psi(n)\sigma(n)}{(p-1)\varphi(n)\big(\psi(n)+\sigma(n)\big)}
-\frac{\psi(n)\big(\varphi(n)+\sigma(n)\big)}{\sigma(n)\big(\varphi(n)+\psi(n)\big)}\nonumber\\
&+\frac{\sigma(n)\big((p-1)\varphi(n)+(p+1)\psi(n)\big)}{\psi(n)\big((p-1)\varphi(n)+(p+1)\sigma(n)\big)}+\frac{\psi(n)\sigma(n)}{(p-1)\varphi(n)\big(\psi(n)+\sigma(n)\big)}
-\frac{\sigma(n)\big(\varphi(n)+\psi(n)\big)}{\psi(n)\big(\varphi(n)+\sigma(n)\big)}\\
&>\frac{9n^2p^2-9np+2}{np(np-1)}\,.
\end{align*}

Case B.  \;  $p\,|\,n$. Using \eqref{phipsisigmanp} and \eqref{phi<psi<sigma} together  lengthy but elementary calculations, we get
\begin{align}\label{CaseB}
&\frac{\big(\sigma(np)-p\sigma(n)\big)\varphi^2(n)}{\psi(n)\big(\varphi(n)+\sigma(n)\big)\big(p\varphi(n)+\sigma(np)\big)}
-\frac{\big(\sigma(np)-p\sigma(n)\big)\psi(n)\varphi(n)}{\sigma(np)\sigma(n)\big(\varphi(n)+\psi(n)\big)}\nonumber\\
&+\frac{\big(\sigma(np)-p\sigma(n)\big)\psi^2(n)}{\varphi(n)\big(\psi(n)+\sigma(n)\big)\big(p\psi(n)+\sigma(np)\big)}\geq0
\end{align}
which, together with \eqref{phipsisigmanp} leads to
\begin{align*}
&\Big(\varphi(np)\big(\psi(np)+\sigma(np)\big)+\psi(np)\big(\varphi(np)+\sigma(np)\big)+\sigma(np)\big(\varphi(np)+\psi(np)\big)\Big)\\
&\times\left(\frac{1}{\varphi(np)\big(\psi(np)+\sigma(np)\big)}+\frac{1}{\psi(np)\big(\varphi(np)+\sigma(np)\big)}+\frac{1}{\sigma(np)\big(\varphi(np)+\psi(np)\big)}\right)\\
&=\Big(p\varphi(n)\big(p\psi(n)+\sigma(np)\big)+p\psi(n)\big(p\varphi(n)+\sigma(np)\big)+p\sigma(np)\big(\varphi(n)+\psi(n)\big)\Big)\\
&\times\left(\frac{1}{p\varphi(n)\big(p\psi(n)+\sigma(np)\big)}+\frac{1}{p\psi(n)\big(p\varphi(n)+\sigma(np)\big)}+\frac{1}{p\sigma(np)\big(\varphi(n)+\psi(n)\big)}\right)\\
&=\bigg(p\Big(\varphi(n)\big(\psi(n)+\sigma(n)\big)+\psi(n)\big(\varphi(n)+\sigma(n)\big)+\sigma(n)\big(\varphi(n)+\psi(n)\big)\Big)\\
&+2\big(\varphi(n)+\psi(n)\big)\big(\sigma(np)-p\sigma(n)\big)\bigg)\frac{1}{p}\\
&\times\left(\frac{1}{\varphi(n)\big(\psi(n)+\sigma(n)\big)}+\frac{1}{\psi(n)\big(\varphi(n)+\sigma(n)\big)}+\frac{1}{\sigma(n)\big(\varphi(n)+\psi(n)\big)}\right)\\
&+\Big(\varphi(n)\big(p\psi(n)+\sigma(np)\big)+\psi(n)\big(p\varphi(n)+\sigma(np)\big)+\sigma(np)\big(\varphi(n)+\psi(n)\big)\Big)\\
&\times\left(\frac{1}{\varphi(n)\big(p\psi(n)+\sigma(np)\big)}+\frac{1}{\psi(n)\big(p\varphi(n)+\sigma(np)\big)}+\frac{1}{\sigma(np)\big(\varphi(n)+\psi(n)\big)}\right)\\
&-\bigg(p\Big(\varphi(n)\big(\psi(n)+\sigma(n)\big)+\psi(n)\big(\varphi(n)+\sigma(n)\big)+\sigma(n)\big(\varphi(n)+\psi(n)\big)\Big)\\
&+2\big(\varphi(n)+\psi(n)\big)\big(\sigma(np)-p\sigma(n)\big)\bigg)\frac{1}{p}\\
&\times\left(\frac{1}{\varphi(n)\big(\psi(n)+\sigma(n)\big)}+\frac{1}{\psi(n)\big(\varphi(n)+\sigma(n)\big)}+\frac{1}{\sigma(n)\big(\varphi(n)+\psi(n)\big)}\right)\\
&\geq\frac{9n^2-9n+2}{n(n-1)}+\frac{2\big(\sigma(np)-p\sigma(n)\big)\varphi^2(n)}{\psi(n)\big(\varphi(n)+\sigma(n)\big)\big(p\varphi(n)+\sigma(np)\big)}\\
&-\frac{2\big(\sigma(np)-p\sigma(n)\big)\psi(n)\varphi(n)}{\sigma(np)\sigma(n)\big(\varphi(n)+\psi(n)\big)}
+\frac{2\big(\sigma(np)-p\sigma(n)\big)\psi^2(n)}{\varphi(n)\big(\psi(n)+\sigma(n)\big)\big(p\psi(n)+\sigma(np)\big)}
>\frac{9n^2p^2-9np+2}{np(np-1)}\,.
\end{align*}
This completes the proof of Theorem \ref{Theorem3}.

\section{Proof of Theorem \ref{Theorem4}}

Consider several cases.

Case 1. \; $\Omega(n)=1$. Bearing in mind that $n$ is a prime number, we have
\begin{align*}
&\frac{\varphi(n)\psi(n)}{\sigma(n)\big(\varphi(n)+\psi(n)\big)}+\frac{\varphi(n)\sigma(n)}{\psi(n)\big(\varphi(n)+\sigma(n)\big)}+\frac{\psi(n)\sigma(n)}{\varphi(n)\big(\psi(n)+\sigma(n)\big)}\\
&=\frac{n-1}{n}+\frac{n+1}{2(n-1)}=\frac{3n^2-3n+2}{2n(n-1)}\,.
\end{align*}

Case 2. \; $\Omega(n)=2$, $n=pq$, where $p$ and $q$ are distinct primes. We write
\begin{align*}
&\frac{\varphi(n)\psi(n)}{\sigma(n)\big(\varphi(n)+\psi(n)\big)}+\frac{\varphi(n)\sigma(n)}{\psi(n)\big(\varphi(n)+\sigma(n)\big)}+\frac{\psi(n)\sigma(n)}{\varphi(n)\big(\psi(n)+\sigma(n)\big)}\\
&=\frac{2(p-1)(q-1)}{(p-1)(q-1)+(p+1)(q+1)}+\frac{(p+1)(q+1)}{2(p-1)(q-1)}\\ 
&=\frac{3p^2q^2-3p^2q-3pq^2+2p^2+10pq+2q^2-3p-3q+3}{2(p-1)(q-1)(pq+1)}\\
&>\frac{3p^2q^2-3pq+2}{2pq(pq-1)}=\frac{3n^2-3n+2}{2n(n-1)}\,.
\end{align*}

Case 3. \; $\Omega(n)=2$, $n=p^2$, where $p$ is a prime. We deduce
\begin{align*}
&\frac{\varphi(n)\psi(n)}{\sigma(n)\big(\varphi(n)+\psi(n)\big)}+\frac{\varphi(n)\sigma(n)}{\psi(n)\big(\varphi(n)+\sigma(n)\big)}+\frac{\psi(n)\sigma(n)}{\varphi(n)\big(\psi(n)+\sigma(n)\big)}\\
&=\frac{p^2(p^2-1)}{(p^2+p+1)(p(p-1)+p(p+1))}+\frac{p(p-1)(p^2+p+1)}{p(p+1)(p(p-1)+p^2+p+1)} \\
&+\frac{p(p+1)(p^2+p+1)}{p(p-1)(p(p+1)+p^2+p+1)}\\
&=\frac{12p^8+24p^7+32p^6+34p^5+37p^4+32p^3+26p^2+14p+5}{8p^8+16p^7+16p^6+4p^5-10p^4-14p^3-12p^2-6p-2}\\
&>\frac{3p^4-3p^2+2}{2p^2(p^2-1)}=\frac{3n^2-3n+2}{2n(n-1)}\,.
\end{align*}
Now we assume that \eqref{lowerbound4} is true for every positive integer $n$ with $\Omega(n)=m$ for some positive integer $m\geq2$.
Let $p$ be a prime number. Then $\Omega(np)=\Omega(n)+1$.

Case A.  \;  $p\nmid n$. Using \eqref{phi<psi<sigma} together with lengthy but elementary calculations, we derive 
\begin{align*}
&\frac{(p-1)\varphi(n)\psi(n)}{\sigma(n)\big((p-1)\varphi(n)+(p+1)\psi(n)\big)}+\frac{\psi(n)\sigma(n)}{(p-1)\varphi(n)\big(\psi(n)+\sigma(n)\big)}
-\frac{\varphi(n)\psi(n)}{\sigma(n)\big(\varphi(n)+\psi(n)\big)}\geq0\,,\\
&\frac{(p-1)\varphi(n)\sigma(n)}{\psi(n)\big((p-1)\varphi(n)+(p+1)\sigma(n)\big)}+\frac{\psi(n)\sigma(n)}{(p-1)\varphi(n)\big(\psi(n)+\sigma(n)\big)}
-\frac{\varphi(n)\sigma(n)}{\psi(n)\big(\varphi(n)+\sigma(n)\big)}\geq0
\end{align*}
which gives us
\begin{align*}
&\frac{\varphi(np)\psi(np)}{\sigma(np)\big(\varphi(np)+\psi(np)\big)}+\frac{\varphi(np)\sigma(np)}{\psi(np)\big(\varphi(np)+\sigma(np)\big)}+\frac{\psi(np)\sigma(np)}{\varphi(np)\big(\psi(np)+\sigma(np)\big)}\\
&=\frac{(p-1)\varphi(n)\psi(n)}{\sigma(n)\big((p-1)\varphi(n)+(p+1)\psi(n)\big)}+\frac{(p-1)\varphi(n)\sigma(n)}{\psi(n)\big((p-1)\varphi(n)+(p+1)\sigma(n)\big)}\\
&+\frac{(p+1)\psi(n)\sigma(n)}{(p-1)\varphi(n)\big(\psi(n)+\sigma(n)\big)}\\
&=\frac{\varphi(n)\psi(n)}{\sigma(n)\big(\varphi(n)+\psi(n)\big)}+\frac{\varphi(n)\sigma(n)}{\psi(n)\big(\varphi(n)+\sigma(n)\big)}+\frac{\psi(n)\sigma(n)}{\varphi(n)\big(\psi(n)+\sigma(n)\big)}\\
&+\frac{(p-1)\varphi(n)\psi(n)}{\sigma(n)\big((p-1)\varphi(n)+(p+1)\psi(n)\big)}+\frac{\psi(n)\sigma(n)}{(p-1)\varphi(n)\big(\psi(n)+\sigma(n)\big)}
-\frac{\varphi(n)\psi(n)}{\sigma(n)\big(\varphi(n)+\psi(n)\big)}\\
&+\frac{(p-1)\varphi(n)\sigma(n)}{\psi(n)\big((p-1)\varphi(n)+(p+1)\sigma(n)\big)}+\frac{\psi(n)\sigma(n)}{(p-1)\varphi(n)\big(\psi(n)+\sigma(n)\big)}
-\frac{\varphi(n)\sigma(n)}{\psi(n)\big(\varphi(n)+\sigma(n)\big)}\\
&\geq\frac{3n^2-3n+2}{2n(n-1)}>\frac{3n^2p^2-3np+2}{2np(np-1)}\,.
\end{align*}

Case B. \;  $p\,|\,n$. From \eqref{CaseB}, we obtain
\begin{align*}
&\frac{\varphi(np)\psi(np)}{\sigma(np)\big(\varphi(np)+\psi(np)\big)}+\frac{\varphi(np)\sigma(np)}{\psi(np)\big(\varphi(np)+\sigma(np)\big)}+\frac{\psi(np)\sigma(np)}{\varphi(np)\big(\psi(np)+\sigma(np)\big)}\\
&=\frac{p\varphi(n)\psi(n)}{\sigma(np)\big(\varphi(n)+\psi(n)\big)}+\frac{\varphi(n)\sigma(np)}{\psi(n)\big(p\varphi(n)+\sigma(np)\big)}+\frac{\psi(n)\sigma(np)}{\varphi(n)\big(p\psi(n)+\sigma(np)\big)}\\
&=\frac{\varphi(n)\psi(n)}{\sigma(n)\big(\varphi(n)+\psi(n)\big)}+\frac{\varphi(n)\sigma(n)}{\psi(n)\big(\varphi(n)+\sigma(n)\big)}+\frac{\psi(n)\sigma(n)}{\varphi(n)\big(\psi(n)+\sigma(n)\big)}\\
&+\frac{\big(\sigma(np)-p\sigma(n)\big)\varphi^2(n)}{\psi(n)\big(\varphi(n)+\sigma(n)\big)\big(p\varphi(n)+\sigma(np)\big)}\\
&-\frac{\big(\sigma(np)-p\sigma(n)\big)\psi(n)\varphi(n)}{\sigma(np)\sigma(n)\big(\varphi(n)+\psi(n)\big)}
+\frac{\big(\sigma(np)-p\sigma(n)\big)\psi^2(n)}{\varphi(n)\big(\psi(n)+\sigma(n)\big)\big(p\psi(n)+\sigma(np)\big)}\\
&\geq\frac{3n^2-3n+2}{2n(n-1)}>\frac{3n^2p^2-3np+2}{2np(np-1)}\,.
\end{align*}
This completes the proof of Theorem \ref{Theorem4}.

\section{Proof of Theorem \ref{Theorem5}}

Consider several cases.

Case 1. \; $\Omega(n)=1$. Taking into account that $n$ is a prime number, we write
\begin{align*}
&\frac{\varphi(n)+\psi(n)}{\varphi(n)+2\sigma(n)+\psi(n)}+\frac{\varphi(n)+\sigma(n)}{\varphi(n)+2\psi(n)+\sigma(n)}+\frac{\psi(n)+\sigma(n)}{\psi(n)+2\varphi(n)+\sigma(n)}\\
&=\frac{2n}{2n+1}+\frac{n+1}{2n}=\frac{6n^2+3n+1}{2n(2n+1)}\,.
\end{align*}

Case 2. \; $\Omega(n)=2$, $n=pq$, where $p$ and $q$ are distinct primes. We have
\begin{align*}
&\frac{\varphi(n)+\psi(n)}{\varphi(n)+2\sigma(n)+\psi(n)}+\frac{\varphi(n)+\sigma(n)}{\varphi(n)+2\psi(n)+\sigma(n)}+\frac{\psi(n)+\sigma(n)}{\psi(n)+2\varphi(n)+\sigma(n)}\\
&=\frac{2\big((p-1)(q-1)+(p+1)(q+1)\big)}{(p-1)(q-1)+3(p+1)(q+1)}+\frac{(p+1)(q+1)}{(p-1)(q-1)+(p+1)(q+1)}\\ 
&=\frac{6p^2q^2+3p^2q+3pq^2+p^2+14pq+q^2+3p+3q+6}{2(pq+1)(2pq+p+q+2)}\\
&>\frac{6p^2q^2+3pq+1}{2pq(2pq+1)}=\frac{6n^2+3n+1}{2n(2n+1)}\,.
\end{align*}

Case 3. \; $\Omega(n)=2$, $n=p^2$, where $p$ is a prime. We get
\begin{align*}
&\frac{\varphi(n)+\psi(n)}{\varphi(n)+2\sigma(n)+\psi(n)}+\frac{\varphi(n)+\sigma(n)}{\varphi(n)+2\psi(n)+\sigma(n)}+\frac{\psi(n)+\sigma(n)}{\psi(n)+2\varphi(n)+\sigma(n)}\\
&=\frac{p(p-1)+p(p+1)}{p(p-1)+2(p^2+p+1)+p(p+1)}+\frac{p(p-1)+p^2+p+1}{p(p-1)+2p(p+1)+p^2+p+1}\\
&+\frac{p(p+1)+p^2+p+1}{p(p+1)+2p(p-1)+p^2+p+1}\\
&=\frac{48p^6+48p^5+68p^4+38p^3+25p^2+6p+2}{32p^6+32p^5+40p^4+20p^3+12p^2+3p+1}\\
&>\frac{6p^4+3p^2+1}{2p^2(2p^2+1)}=\frac{6n^2+3n+1}{2n(2n+1)}\,.
\end{align*}
Now we assume that \eqref{lowerbound5} is true for every positive integer $n$ with $\Omega(n)=m$ for some positive integer $m\geq2$.
Let $p$ be a prime number. Then $\Omega(np)=\Omega(n)+1$.

Case A.  \;  $p\nmid n$. It follows easily from \eqref{phi<psi<sigma} that
\begin{align*}
&\frac{1}{\big((p+1)\psi(n)+2(p-1)\varphi(n)+(p+1)\sigma(n)\big)\big(\psi(n)+2\varphi(n)+\sigma(n)\big)}\\
&-\frac{1}{\big((p-1)\varphi(n)+2(p+1)\sigma(n)+(p+1)\psi(n)\big)\big(\varphi(n)+2\sigma(n)+\psi(n)\big)}\geq0\,,\\
&\frac{1}{\big((p+1)\psi(n)+2(p-1)\varphi(n)+(p+1)\sigma(n)\big)\big(\psi(n)+2\varphi(n)+\sigma(n)\big)}\\
&-\frac{1}{\big((p-1)\varphi(n)+2(p+1)\psi(n)+(p+1)\sigma(n)\big)\big(\varphi(n)+2\psi(n)+\sigma(n)\big)}\geq0
\end{align*}
which  implies
\begin{align*}
&\frac{\varphi(np)+\psi(np)}{\varphi(np)+2\sigma(np)+\psi(np)}+\frac{\varphi(np)+\sigma(np)}{\varphi(np)+2\psi(np)+\sigma(np)}+\frac{\psi(np)+\sigma(np)}{\psi(np)+2\varphi(np)+\sigma(np)}\\
&=\frac{(p-1)\varphi(n)+(p+1)\psi(n)}{(p-1)\varphi(n)+2(p+1)\sigma(n)+(p+1)\psi(n)}+\frac{(p-1)\varphi(n)+(p+1)\sigma(n)}{(p-1)\varphi(n)+2(p+1)\psi(n)+(p+1)\sigma(n)}\\
&+\frac{(p+1)\big(\psi(n)+\sigma(n)\big)}{(p+1)\psi(n)+2(p-1)\varphi(n)+(p+1)\sigma(n)}\\
&=\frac{\varphi(n)+\psi(n)}{\varphi(n)+2\sigma(n)+\psi(n)}+\frac{\varphi(n)+\sigma(n)}{\varphi(n)+2\psi(n)+\sigma(n)}+\frac{\psi(n)+\sigma(n)}{\psi(n)+2\varphi(n)+\sigma(n)}\\
&+\frac{(p-1)\varphi(n)+(p+1)\psi(n)}{(p-1)\varphi(n)+2(p+1)\sigma(n)+(p+1)\psi(n)}-\frac{\varphi(n)+\psi(n)}{\varphi(n)+2\sigma(n)+\psi(n)}\\
&+\frac{(p-1)\varphi(n)+(p+1)\sigma(n)}{(p-1)\varphi(n)+2(p+1)\psi(n)+(p+1)\sigma(n)}-\frac{\varphi(n)+\sigma(n)}{\varphi(n)+2\psi(n)+\sigma(n)}\\
&+\frac{(p+1)\big(\psi(n)+\sigma(n)\big)}{(p+1)\psi(n)+2(p-1)\varphi(n)+(p+1)\sigma(n)}-\frac{\psi(n)+\sigma(n)}{\psi(n)+2\varphi(n)+\sigma(n)}\\
&\geq\frac{6n^2+3n+1}{2n(2n+1)}+\frac{4\varphi(n)\big(\psi(n)+\sigma(n)\big)}{\big((p+1)\psi(n)+2(p-1)\varphi(n)+(p+1)\sigma(n)\big)\big(\psi(n)+2\varphi(n)+\sigma(n)\big)}\\
&-\frac{4\varphi(n)\sigma(n)}{\big((p-1)\varphi(n)+2(p+1)\sigma(n)+(p+1)\psi(n)\big)\big(\varphi(n)+2\sigma(n)+\psi(n)\big)}\\
&-\frac{4\varphi(n)\psi(n)}{\big((p-1)\varphi(n)+2(p+1)\psi(n)+(p+1)\sigma(n)\big)\big(\varphi(n)+2\psi(n)+\sigma(n)\big)}\\
&>\frac{6n^2p^2+3np+1}{2np(2np+1)}\,.
\end{align*}

Case B.  \;  $p\,|\,n$. By \eqref{phipsisigmanp} and \eqref{phi<psi<sigma}, we deduce
\begin{align*}
&\frac{1}{\big(p\varphi(n)+2p\psi(n)+\sigma(np)\big)\big(\varphi(n)+2\psi(n)+\sigma(n)\big)}\\
&-\frac{1}{\big(p\varphi(n)+2\sigma(np)+p\psi(n)\big)\big(\varphi(n)+2\sigma(n)+\psi(n)\big)}\geq0\,,\\
&\frac{1}{\big(p\psi(n)+2p\varphi(n)+\sigma(np)\big)\big(\psi(n)+2\varphi(n)+\sigma(n)\big)}\\
&-\frac{1}{\big(p\varphi(n)+2\sigma(np)+p\psi(n)\big)\big(\varphi(n)+2\sigma(n)+\psi(n)\big)}\geq0
\end{align*}
which, together with \eqref{phipsisigmanp} yields
\begin{align*}
&\frac{\varphi(np)+\psi(np)}{\varphi(np)+2\sigma(np)+\psi(np)}+\frac{\varphi(np)+\sigma(np)}{\varphi(np)+2\psi(np)+\sigma(np)}+\frac{\psi(np)+\sigma(np)}{\psi(np)+2\varphi(np)+\sigma(np)}\\
\end{align*}
\begin{align*}
&=\frac{p\big(\varphi(n)+\psi(n)\big)}{p\varphi(n)+2\sigma(np)+p\psi(n)}+\frac{p\varphi(n)+\sigma(np)}{p\varphi(n)+2p\psi(n)+\sigma(np)}+\frac{p\psi(n)+\sigma(np)}{p\psi(n)+2p\varphi(n)+\sigma(np)}\\
&=\frac{\varphi(n)+\psi(n)}{\varphi(n)+2\sigma(n)+\psi(n)}+\frac{\varphi(n)+\sigma(n)}{\varphi(n)+2\psi(n)+\sigma(n)}+\frac{\psi(n)+\sigma(n)}{\psi(n)+2\varphi(n)+\sigma(n)}\\
&+\frac{p\big(\varphi(n)+\psi(n)\big)}{p\varphi(n)+2\sigma(np)+p\psi(n)}-\frac{\varphi(n)+\psi(n)}{\varphi(n)+2\sigma(n)+\psi(n)}\\
&+\frac{p\varphi(n)+\sigma(np)}{p\varphi(n)+2p\psi(n)+\sigma(np)}-\frac{\varphi(n)+\sigma(n)}{\varphi(n)+2\psi(n)+\sigma(n)}\\
&+\frac{p\psi(n)+\sigma(np)}{p\psi(n)+2p\varphi(n)+\sigma(np)}-\frac{\psi(n)+\sigma(n)}{\psi(n)+2\varphi(n)+\sigma(n)}\\
&\geq\frac{6n^2+3n+1}{2n(2n+1)}+\frac{2\big(p\sigma-\sigma(np)\big)\big(\varphi(n)+\psi(n)\big)}{\big(p\varphi(n)+2\sigma(np)+p\psi(n)\big)\big(\varphi(n)+2\sigma(n)+\psi(n)\big)}\\
&+\frac{2\big(\sigma(np)-p\sigma\big)\psi(n)}{\big(p\varphi(n)+2p\psi(n)+\sigma(np)\big)\big(\varphi(n)+2\psi(n)+\sigma(n)\big)}\\
&+\frac{2\big(\sigma(np)-p\sigma\big)\varphi(n)}{\big(p\psi(n)+2p\varphi(n)+\sigma(np)\big)\big(\psi(n)+2\varphi(n)+\sigma(n)\big)}\\
&>\frac{6n^2p^2+3np+1}{2np(2np+1)}\,.
\end{align*}
This completes the proof of Theorem \ref{Theorem5}.

\section{Proof of Theorem \ref{Theorem6}}

Consider several cases.

Case 1. \; $\Omega(n)=1$. Taking into account that $n$ is a prime number, we have
\begin{align*}
&\frac{\varphi(n)\big(\psi(n)+\sigma(n)\big)}{\psi(n)\big(\varphi(n)+\sigma(n)\big)+\sigma(n)\big(\varphi(n)+\psi(n)\big)}
+\frac{\psi(n)\big(\varphi(n)+\sigma(n)\big)}{\varphi(n)\big(\psi(n)+\sigma(n)\big)+\sigma(n)\big(\varphi(n)+\psi(n)\big)}\\
&+\frac{\sigma(n)\big(\varphi(n)+\psi(n)\big)}{\varphi(n)\big(\psi(n)+\sigma(n)\big)+\psi(n)\big(\varphi(n)+\sigma(n)\big)}=\frac{n-1}{2n}+\frac{2n}{2n-1}=\frac{6n^2-3n+1}{2n(2n-1)}\,.
\end{align*}

Case 2. \; $\Omega(n)=2$, $n=pq$, where $p$ and $q$ are distinct primes. We write
\begin{align*}
&\frac{\varphi(n)\big(\psi(n)+\sigma(n)\big)}{\psi(n)\big(\varphi(n)+\sigma(n)\big)+\sigma(n)\big(\varphi(n)+\psi(n)\big)}
+\frac{\psi(n)\big(\varphi(n)+\sigma(n)\big)}{\varphi(n)\big(\psi(n)+\sigma(n)\big)+\sigma(n)\big(\varphi(n)+\psi(n)\big)}\\
&+\frac{\sigma(n)\big(\varphi(n)+\psi(n)\big)}{\varphi(n)\big(\psi(n)+\sigma(n)\big)+\psi(n)\big(\varphi(n)+\sigma(n)\big)}\\
\end{align*}
\begin{align*}
&=\frac{(p-1)(q-1)}{(p-1)(q-1)+(p+1)(q+1)}+\frac{2\big((p-1)(q-1)+(p+1)(q+1)\big)}{3(p-1)(q-1)+(p+1)(q+1)}\\ 
&=\frac{6p^2q^2-3p^2q-3pq^2+p^2+14pq+q^2-3p-3q+6}{2(pq+1)(2pq-p-q+2)}
>\frac{6p^2q^2-3pq+1}{2pq(2pq-1)}=\frac{6n^2-3n+1}{2n(2n-1)}\,.
\end{align*}

Case 3. \; $\Omega(n)=2$, $n=p^2$, where $p$ is a prime. We derive
\begin{align*}
&\frac{\varphi(n)\big(\psi(n)+\sigma(n)\big)}{\psi(n)\big(\varphi(n)+\sigma(n)\big)+\sigma(n)\big(\varphi(n)+\psi(n)\big)}
+\frac{\psi(n)\big(\varphi(n)+\sigma(n)\big)}{\varphi(n)\big(\psi(n)+\sigma(n)\big)+\sigma(n)\big(\varphi(n)+\psi(n)\big)}\\
&+\frac{\sigma(n)\big(\varphi(n)+\psi(n)\big)}{\varphi(n)\big(\psi(n)+\sigma(n)\big)+\psi(n)\big(\varphi(n)+\sigma(n)\big)}\\
&=\frac{p(p-1)\big(p(p+1)+p^2+p+1\big)}{p(p+1)\big(p(p-1)+p^2+p+1\big)+(p^2+p+1)\big(p(p-1)+p(p+1)\big)}+\\
&+\frac{p(p+1)\big(p(p-1)+p^2+p+1\big)}{p(p-1)\big(p(p+1)+p^2+p+1\big)+(p^2+p+1)\big(p(p-1)+p(p+1)\big)}+\\
&+\frac{(p^2+p+1)\big(p(p-1)+p(p+1)\big)}{p(p-1)\big(p(p+1)+p^2+p+1\big)+p(p+1)\big(p(p-1)+p^2+p+1\big)}\\
&=\frac{24p^7+36p^6+40p^5+25p^4+15p^3+4p^2+p-1}{16p^7+24p^6+24p^5+10p^4+p^3-2p^2-p}\\
&>\frac{6p^4-3p^2+1}{2p^2(2p^2-1)}=\frac{6n^2-3n+1}{2n(2n-1)}\,.
\end{align*}
Now we assume that \eqref{lowerbound6} is true for every positive integer $n$ with $\Omega(n)=m$ for some positive integer $m\geq2$.
Let $p$ be a prime number. Then $\Omega(np)=\Omega(n)+1$.

Case A.  \;  $p\nmid n$. Put

\begin{align*}
\Delta_1&=\Big(\psi(n)\big((p-1)\varphi(n)+(p+1)\sigma(n)\big)+\sigma(n)\big((p-1)\varphi(n)+(p+1)\psi(n)\big)\Big)\\
&\times\Big(\psi(n)\big(\varphi(n)+\sigma(n)\big)+\sigma(n)\big(\varphi(n)+\psi(n)\big)\Big)\,,\\
\Delta_2&=\Big((p-1)\varphi(n)\big(\psi(n)+\sigma(n)\big)+\sigma(n)\big((p-1)\varphi(n)+(p+1)\psi(n)\big)\Big)\\
&\times\Big(\varphi(n)\big(\psi(n)+\sigma(n)\big)+\sigma(n)\big(\varphi(n)+\psi(n)\big)\Big)\,,\\
\Delta_3&=\Big((p-1)\varphi(n)\big(\psi(n)+\sigma(n)\big)+\psi(n)\big((p-1)\varphi(n)+(p+1)\sigma(n)\big)\Big)\\
&\times\Big(\varphi(n)\big(\psi(n)+\sigma(n)\big)+\psi(n)\big(\varphi(n)+\sigma(n)\big)\Big)\,.
\end{align*}
Using \eqref{phi<psi<sigma}, we obtain  
\begin{equation*}
\Delta_1>\Delta_2\,, \qquad \Delta_1>\Delta_3
\end{equation*}
which leads to
\begin{align*}
&\frac{\varphi(np)\big(\psi(np)+\sigma(np)\big)}{\psi(np)\big(\varphi(np)+\sigma(np)\big)+\sigma(np)\big(\varphi(np)+\psi(np)\big)}\\
&+\frac{\psi(np)\big(\varphi(np)+\sigma(np)\big)}{\varphi(np)\big(\psi(np)+\sigma(np)\big)+\sigma(np)\big(\varphi(np)+\psi(np)\big)}\\
&+\frac{\sigma(np)\big(\varphi(np)+\psi(np)\big)}{\varphi(np)\big(\psi(np)+\sigma(np)\big)+\psi(np)\big(\varphi(np)+\sigma(np)\big)}\\
&=\frac{(p-1)\varphi(n)\big(\psi(n)+\sigma(n)\big)}{\psi(n)\big((p-1)\varphi(n)+(p+1)\sigma(n)\big)+\sigma(n)\big((p-1)\varphi(n)+(p+1)\psi(n)\big)}\\
&+\frac{\psi(n)\big((p-1)\varphi(n)+(p+1)\sigma(n)\big)}{(p-1)\varphi(n)\big(\psi(n)+\sigma(n)\big)+\sigma(n)\big((p-1)\varphi(n)+(p+1)\psi(n)\big)}\\
&+\frac{\sigma(n)\big((p-1)\varphi(n)+(p+1)\psi(n)\big)}{(p-1)\varphi(n)\big(\psi(n)+\sigma(n)\big)+\psi(n)\big((p-1)\varphi(n)+(p+1)\sigma(n)\big)}\\
&=\frac{\varphi(n)\big(\psi(n)+\sigma(n)\big)}{\psi(n)\big(\varphi(n)+\sigma(n)\big)+\sigma(n)\big(\varphi(n)+\psi(n)\big)}
+\frac{\psi(n)\big(\varphi(n)+\sigma(n)\big)}{\varphi(n)\big(\psi(n)+\sigma(n)\big)+\sigma(n)\big(\varphi(n)+\psi(n)\big)}\\
&+\frac{\sigma(n)\big(\varphi(n)+\psi(n)\big)}{\varphi(n)\big(\psi(n)+\sigma(n)\big)+\psi(n)\big(\varphi(n)+\sigma(n)\big)}\\
&+\frac{(p-1)\varphi(n)\big(\psi(n)+\sigma(n)\big)}{\psi(n)\big((p-1)\varphi(n)+(p+1)\sigma(n)\big)+\sigma(n)\big((p-1)\varphi(n)+(p+1)\psi(n)\big)}\\
&-\frac{\varphi(n)\big(\psi(n)+\sigma(n)\big)}{\psi(n)\big(\varphi(n)+\sigma(n)\big)+\sigma(n)\big(\varphi(n)+\psi(n)\big)}\\
&+\frac{\psi(n)\big((p-1)\varphi(n)+(p+1)\sigma(n)\big)}{(p-1)\varphi(n)\big(\psi(n)+\sigma(n)\big)+\sigma(n)\big((p-1)\varphi(n)+(p+1)\psi(n)\big)}\\
&-\frac{\psi(n)\big(\varphi(n)+\sigma(n)\big)}{\varphi(n)\big(\psi(n)+\sigma(n)\big)+\sigma(n)\big(\varphi(n)+\psi(n)\big)}\\
&+\frac{\sigma(n)\big((p-1)\varphi(n)+(p+1)\psi(n)\big)}{(p-1)\varphi(n)\big(\psi(n)+\sigma(n)\big)+\psi(n)\big((p-1)\varphi(n)+(p+1)\sigma(n)\big)}\\
&-\frac{\sigma(n)\big(\varphi(n)+\psi(n)\big)}{\varphi(n)\big(\psi(n)+\sigma(n)\big)+\psi(n)\big(\varphi(n)+\sigma(n)\big)}\\
&\geq\frac{6n^2-3n+1}{2n(2n-1)}-\frac{4\sigma(n)\psi(n)\varphi(n)\big(\psi(n)+\sigma(n)\big)}{\Delta_1}+\frac{4\sigma^2(n)\psi(n)\varphi(n)}{\Delta_2}+\frac{4\sigma(n)\psi^2(n)\varphi(n)}{\Delta_3}\\
&>\frac{6n^2p^2-3np+1}{2np(2np-1)}-\frac{4\sigma(n)\psi(n)\varphi(n)\big(\psi(n)+\sigma(n)\big)}{\Delta_1}+\frac{4\sigma^2(n)\psi(n)\varphi(n)}{\Delta_1}+\frac{4\sigma(n)\psi^2(n)\varphi(n)}{\Delta_1}\\
&=\frac{6n^2p^2-3np+1}{2np(2np-1)}\,.
\end{align*}

Case B. \;  $p\,|\,n$. Set
\begin{align*}
\Omega_1&=\Big(\psi(n)\big(p\varphi(n)+\sigma(np)\big)+\sigma(np)\big(\varphi(n)+\psi(n)\big)\Big)\\
&\times\Big(\psi(n)\big(\varphi(n)+\sigma(n)\big)+\sigma(n)\big(\varphi(n)+\psi(n)\big)\Big)\,,\\
\Omega_2&=\Big(\varphi(n)\big(p\psi(n)+\sigma(np)\big)+\sigma(np)\big(\varphi(n)+\psi(n)\big)\Big)\\
&\times\Big(\varphi(n)\big(\psi(n)+\sigma(n)\big)+\sigma(n)\big(\varphi(n)+\psi(n)\big)\Big)\,,\\
\Omega_3&=\Big(\varphi(n)\big(p\psi(n)+\sigma(np)\big)+\psi(n)\big(p\varphi(n)+\sigma(np)\big)\Big)\\
&\times\Big(\varphi(n)\big(\psi(n)+\sigma(n)\big)+\psi(n)\big(\varphi(n)+\sigma(n)\big)\Big)\,.
\end{align*}
From \eqref{phipsisigmanp} and \eqref{phi<psi<sigma}, we get
\begin{equation*}
\Omega_1>\Omega_3\,, \qquad \Omega_2>\Omega_3
\end{equation*}
which, together with \eqref{phipsisigmanp} gives us
\begin{align*}
&\frac{\varphi(np)\big(\psi(np)+\sigma(np)\big)}{\psi(np)\big(\varphi(np)+\sigma(np)\big)+\sigma(np)\big(\varphi(np)+\psi(np)\big)}\\
&+\frac{\psi(np)\big(\varphi(np)+\sigma(np)\big)}{\varphi(np)\big(\psi(np)+\sigma(np)\big)+\sigma(np)\big(\varphi(np)+\psi(np)\big)}\\
&+\frac{\sigma(np)\big(\varphi(np)+\psi(np)\big)}{\varphi(np)\big(\psi(np)+\sigma(np)\big)+\psi(np)\big(\varphi(np)+\sigma(np)\big)}\\
&=\frac{\varphi(n)\big(p\psi(n)+\sigma(np)\big)}{\psi(n)\big(p\varphi(n)+\sigma(np)\big)+\sigma(np)\big(\varphi(n)+\psi(n)\big)}\\
&+\frac{\psi(n)\big(p\varphi(n)+\sigma(np)\big)}{\varphi(n)\big(p\psi(n)+\sigma(np)\big)+\sigma(np)\big(\varphi(n)+\psi(n)\big)}\\
&+\frac{\sigma(np)\big(\varphi(n)+\psi(n)\big)}{\varphi(n)\big(p\psi(n)+\sigma(np)\big)+\psi(n)\big(p\varphi(n)+\sigma(np)\big)}\\
&=\frac{\varphi(n)\big(\psi(n)+\sigma(n)\big)}{\psi(n)\big(\varphi(n)+\sigma(n)\big)+\sigma(n)\big(\varphi(n)+\psi(n)\big)}
+\frac{\psi(n)\big(\varphi(n)+\sigma(n)\big)}{\varphi(n)\big(\psi(n)+\sigma(n)\big)+\sigma(n)\big(\varphi(n)+\psi(n)\big)}\\
&+\frac{\sigma(n)\big(\varphi(n)+\psi(n)\big)}{\varphi(n)\big(\psi(n)+\sigma(n)\big)+\psi(n)\big(\varphi(n)+\sigma(n)\big)}\\
&+\frac{\varphi(n)\big(p\psi(n)+\sigma(np)\big)}{\psi(n)\big(p\varphi(n)+\sigma(np)\big)+\sigma(np)\big(\varphi(n)+\psi(n)\big)}
-\frac{\varphi(n)\big(\psi(n)+\sigma(n)\big)}{\psi(n)\big(\varphi(n)+\sigma(n)\big)+\sigma(n)\big(\varphi(n)+\psi(n)\big)}\\
&+\frac{\psi(n)\big(p\varphi(n)+\sigma(np)\big)}{\varphi(n)\big(p\psi(n)+\sigma(np)\big)+\sigma(np)\big(\varphi(n)+\psi(n)\big)}
-\frac{\psi(n)\big(\varphi(n)+\sigma(n)\big)}{\varphi(n)\big(\psi(n)+\sigma(n)\big)+\sigma(n)\big(\varphi(n)+\psi(n)\big)}\\
\end{align*}
\begin{align*}
&+\frac{\sigma(np)\big(\varphi(n)+\psi(n)\big)}{\varphi(n)\big(p\psi(n)+\sigma(np)\big)+\psi(n)\big(p\varphi(n)+\sigma(np)\big)}
-\frac{\sigma(n)\big(\varphi(n)+\psi(n)\big)}{\varphi(n)\big(\psi(n)+\sigma(n)\big)+\psi(n)\big(\varphi(n)+\sigma(n)\big)}\\
&\geq\frac{6n^2-3n+1}{2n(2n-1)}-\frac{2\varphi(n)\psi^2(n)\big(\sigma(np)-p\sigma\big)}{\Omega_1}-\frac{2\varphi^2(n)\psi(n)\big(\sigma(np)-p\sigma\big)}{\Omega_2}\\
&+\frac{2\big(\varphi(n)\psi^2(n)+\varphi^2(n)\psi(n)\big)\big(\sigma(np)-p\sigma\big)}{\Omega_3}\\
&>\frac{6n^2p^2-3np+1}{2np(2np-1)}-\frac{2\varphi(n)\psi^2(n)\big(\sigma(np)-p\sigma\big)}{\Omega_1}-\frac{2\varphi^2(n)\psi(n)\big(\sigma(np)-p\sigma\big)}{\Omega_2}\\
&+\frac{2\varphi(n)\psi^2(n)\big(\sigma(np)-p\sigma\big)}{\Omega_1}+\frac{2\varphi^2(n)\psi(n)\big(\sigma(np)-p\sigma\big)}{\Omega_2}\\
&=\frac{6n^2p^2-3np+1}{2np(2np-1)}\,.
\end{align*}
This completes the proof of Theorem \ref{Theorem6}.

\vskip20pt
\footnotesize
\begin{flushleft}
S. I. Dimitrov\\
\quad\\
Faculty of Applied Mathematics and Informatics\\
Technical University of Sofia \\
Blvd. St. Kliment Ohridski 8 \\
Sofia 1000, Bulgaria\\
e-mail: sdimitrov@tu-sofia.bg\\
\end{flushleft}

\end{document}